\documentclass[reqno, 11pt]{amsart}
\usepackage{amssymb, mathtools}
\usepackage{mathrsfs}
\usepackage{euscript}
\usepackage{graphicx}
\usepackage[left]{lineno}
\usepackage{xcolor}
\usepackage{float}
\usepackage{pgfplots, caption, subcaption, tikz, tikz-cd} % packages for pictures

\usetikzlibrary{calc,decorations.pathmorphing}

\usepackage[T1]{fontenc}
\usepackage[utf8]{inputenc}

\newtheorem{theorem}{Theorem}[section]

\newtheorem{proposition}[theorem]{Proposition}

\newtheorem{definition}[theorem]{Definition}

\theoremstyle{remark}
\newtheorem{remark}[theorem]{Remark}

\theoremstyle{assumption}

\numberwithin{equation}{section}

\DeclareMathOperator*{\Lip}{\mathrm{Lip}}

\title[Absolutely minimal extensions in finite metric spaces]
      {A note on absolutely minimal extensions \\ in finite metric spaces}\thanks{The authors were partially supported by the Austrian Science Foundation (FWF) under grants P 36344-N and F 100800.}

\author[Alberto Dom\'inguez Corella]{Alberto Dom\'inguez Corella}
\address[Alberto Dom\'inguez Corella]{Institut f\"ur Mathematik und Wissenschaftliches Rechnen,
Universit\"at Graz, Heinrichstrasse 36, 8010 Graz, Austria}
\address{Institut f\"ur Stochastik und Wirtschaftsmathematik,
Technische Universit\"at Wien,
Wiedner Hauptstrasse 8, VADOR E105-04, 1040 Wien, Austria}

\email{alberto.of.sonora@gmail.com}
\urladdr{https://sites.google.com/view/alberto-of-sonora}

\author[Tr\'i Minh L\^e]{Tr\'i Minh L\^e}
\address[Tr\'i Minh L\^e]{Institut f\"ur Stochastik und Wirtschaftsmathematik,
Technische Universit\"at Wien,
Wiedner Hauptstrasse 8, VADOR E105-04, 1040 Wien, Austria}
\email{minh.le@tuwien.ac.at}
\urladdr{https://sites.google.com/view/tri-minh-le}

\keywords{Finite metric spaces, Absolutely minimal Lipschitz extensions, McShane--Whitney formulas}
\subjclass[2020]{26A16, 39B82, 05C12}

\begin{document}
\maketitle
\begin{abstract}\vspace{-3em}
Absolutely minimal Lipschitz extensions (AMLEs) are known to exist in many infinite metric settings, but the finite case is less settled. In metric spaces with at most four points, every function on a nonempty subset admits an AMLE in the sense that the Lipschitz constant cannot be further reduced on sets that are disjoint from the prescribed domain. We show that in five-point spaces such extensions may fail to exist.
\iffalse
Absolutely minimal Lipschitz extensions (AMLEs)  are known to exist in many infinite metric settings, but the finite case is less settled.
In metric spaces with at most four points, every function on a nonempty subset admits an AMLE in the sense that the Lipschitz constant cannot be further reduced on sets disjoint from the prescribed domain. 
We show that in five–point spaces such extensions may fail to exist.
\fi 
\end{abstract}
\section{Introduction, existence failure and state of the art}

It is well known that in any metric space, a Lipschitz function on a subset can be extended to the entire space while preserving its Lipschitz constant.
Such extensions are referred to as minimal Lipschitz extensions, as they attain the smallest possible Lipschitz constant among all Lipschitz extensions, see \cite{McShane, Whitney}.

\medskip

More precisely, let $(X, d)$ be a metric space, let $A \subset X$ be a closed subset and let $f: A \to \mathbb{R}$ be a Lipschitz function on $A$.
A function $g:X\to\mathbb R$ is a \textit{minimal Lipschitz extension} if $g|_A = f$ and $\Lip\big(f, A \big) = \Lip\big(g, X \big)$,  where $\Lip(f, A)$ denotes the Lipschitz constant of $f$ on $A$:
\[
	\Lip(f, A) := \sup_{\substack{x, y \in A \\ x \neq y} } \dfrac{|f(x) - f(y)|}{d(x, y)}.
\]
Two natural candidates for such extensions are the so--called McShane--Whitney extensions
\begin{equation}\label{mM}
		m(x):=\sup_{y\in A}\{f(y) - \Lip(f, A) d(y,x)\}\quad \text{and}\quad  M(x):=\inf_{y\in A}\{f(y) + \Lip(f, A) d(y,x)\}.
\end{equation}

\noindent

The multiplicity and instability of these extensions are somewhat unsatisfactory; additionally the definition involves only the global Lipschitz constant of the extension and ignores what may happen locally. 
For this reason, the definition of \textit{absolutely minimal Lipschitz extensions} was introduced in Euclidean spaces by Aronsson \cite{A_65, A_66, A_67}:

\medskip

\textit{Let $\Omega \subset \mathbb{R}^n$ be a nonempty compact set and let $f: \partial \Omega \to \mathbb{R}$ be a Lipschitz function. Then there exists a minimal Lipschitz extension $u : \Omega \to \mathbb{R}$ of $f$ such that  
\begin{equation}\label{amle}
	\Lip(u, V) = \Lip(u, \partial V) \quad \text{for every open set } V \subset\subset \Omega.
\end{equation}  
Such a function $u$ is called an AMLE of $f$.  
}

\medskip

It has been shown that AMLEs satisfy the so-called infinity Laplace equation and are strongly connected to tug-of-war games, see e.g.,   \cite{ACJ_2004, BCR_2023, C_2019, CP_2007, Gruyer, PSSW_2009}.  
This definition, however, does not extend naturally to general metric spaces, where boundaries may be empty. 
Still, many efforts have been made to develop suitable notions of AMLEs in the metric setting, for instance in complete length spaces and in metric measure spaces, see \cite{DLV_2024, DJS_2019, JS_2006}.

\medskip

To tackle this difficulty, Juutinen \cite{J_2002} introduced an alternative definition of AMLE that does not rely on the boundary. The central idea is that extensions be locally Lipschitz minimal with respect to suitable test functions.

	\begin{definition}\label{def.amle.juuti}
	Let $(X, d)$ be a metric space, let $A$ be a nonempty closed subset of $X$ and let $f: A \to \mathbb R$ be a Lipschitz function.
	We say that $h:X\to\mathbb R$ is an absolutely minimal Lipschitz extension (AMLE)  relative to a family $\EuScript{U} \subset 2^X$ if
	\begin{itemize}
		\item[$(i)$] $h$ is a minimal Lipschitz extension of $f$.
		
		\item[$(ii)$] For any $U \in \EuScript U$ and any minimal extension  $g:X\to\mathbb R$ of $f$ satisfying  $g|_{X\setminus U} = h|_{X\setminus U}$, there holds 
		\[
		\Lip\big(h, U\big)\le \Lip\big(g,U\big).
		\] 
	\end{itemize}
\end{definition}

Different choices of $\EuScript U$ give rise to different variants of AMLE:
\begin{itemize}
  \item If $\EuScript U = 2^X$ is the family of all subsets of $X$, one recovers Juutinen’s definition of AMLE \cite{J_2002}.
  \item If $\EuScript U = 2^{X \setminus A}$ is the family of subsets contained in $X \setminus A$, one obtains the relaxed formulation considered by Mazón, Rossi and Toledo \cite{MRT_2012} in order to study AMLE in particular discrete spaces.
  \item If $X$ is a complete length space and $\EuScript U \subset 2^{X \setminus A}$ is the family of open subsets whose closures are contained in $X \setminus A$, then Definition~\ref{def.amle.juuti}--\em{(ii)} is equivalent to: $\Lip(h, U) = \Lip(h, \partial U)$ for every $U \in \EuScript U$, see \cite[Proposition 2.6]{DLV_2024}.
\end{itemize}
 
Juutinen proved that in separable length spaces AMLEs relative to $2^X$ always exist. However, he made no comment on denumerable metric spaces (finite or countable), apart from giving a trivial three-point example in which extension is not possible \cite[Example 2.2]{J_2002}. 
He also noted there that one might relax the definition to AMLE relative to $2^{X \setminus A}$ to ensure existence in his example.  
This is retaken by Mazón, Rossi and Toledo in \cite[Definition 2.1]{MRT_2012} to address problems in a discrete setting.

\medskip
 
In Juutinen’s three–point example, AMLEs relative to $2^{X \setminus A}$ do exist, although they are not unique. 
On the other hand, for four–point metric spaces one can still verify existence of such extensions. For completeness, we include a short proof of this fact, as this does not appear in the current literature and give an example of non-uniqueness in Remark \ref{rmk01}.

 	\begin{proposition}\label{p0}
 		If $X$ consists of four points and $A \subset X$ is nonempty, then any function $f:A \to \mathbb R$ admits an AMLE relative to $2^{X \setminus A}$.
 	\end{proposition}
Existence in this case is encouraging in finite settings;  however after more than twenty years, few comments in this regard have appeared, in particular it has not been addressed whether there exist AMLEs relative to $2^{X \setminus A}$ in finite metric spaces. 

\medskip

The main result of this note is a construction showing that AMLEs relative to $2^{X \setminus A}$ need not exist when metric spaces have more than four points.
 	\begin{theorem}\label{thm0}
 	There is a five–point metric space with  a function defined on a two-point subset for which no AMLE relative to $2^{X \setminus A}$ exists. 
 	\end{theorem}

We conclude with some remarks.
\begin{remark} $\, $
\begin{enumerate}
    \item[\textit{(i)}]  The constructed metric space in the proof of Theorem \ref{thm0} can be embedded into $\mathbb R^d$ for every $d \geq 1$. We also give a non-Euclidean construction in Remark \ref{rdms}.
    \smallbreak
    \item[\textit{(ii)}] Adapting the proof of Theorem~\ref{thm0}, one can construct counterexamples in arbitrary finite metric spaces with more than five points. 

     \item[\textit{(iii)}] The construction given in the proof of Theorem~\ref{thm0} also shows that AMLE's in the sense of Juutinen \cite{J_2002} (relative to $2^{X}$) also fail to exist in discrete settings.
\end{enumerate}

\end{remark}

\section{Proofs}
\subsection{Proof of Proposition \ref{p0}}
If $A$ has one, three, or four points, the result is trivial; therefore we assume $A$ has two elements.
Consider $X = \{ x_1, \, x_2, \, x_3, \, x_4 \}$ and $A = \{ x_1, \, x_2 \}$. 
Fix a function $f:  A \to \mathbb{R}$.
Define 
	\begin{align}\label{minseq}
	    \mathcal{L}(\{ x_3, \, x_4 \}) = \inf \left\{ \Lip(g, \{x_3, \, x_4\}): \, \text{$g$ is a minimal Lipschitz extension of $f$} \right\}.
	\end{align}
Let \((g_n)_{n\in\mathbb N}\) be a minimizing sequence for (\ref{minseq}), i.e., each $g_n:X\to\mathbb R$ is a minimal Lipschitz extension of \(f\) and 
\[
\operatorname{Lip}(g_n,\{x_3, \, x_4\}) \longrightarrow \mathcal{L}(\{x_3, \, x_4\})
\quad\text{as } n \longrightarrow +\infty.
\]
Due to the  McShane--Whitney bounds, the set over which the minimization occurs is bounded; so  up to a subsequence, there exists $\bar g$  such that $g_n \to \bar g$   as $n \to \infty$.
It is straightforward to see that $\mathcal{L}(\{ x_3, \, x_4 \}) = \Lip(\bar g, \{ x_3, \, x_4\})$ and that $\bar{g}$ is a minimal Lipschitz extension of $f$.
We conclude that $\bar g$ is an AMLE of $f$ relative to $2^{X \setminus A}$.\qed

\begin{remark}
  Using the same reasoning as in the proof of Proposition~\ref{p0}, it can be seen that if 
\(X\) is an \(n\)-point metric space and \(A \subset X\) has \(n-2\) points, then every 
function \(f:A \to \mathbb{R}\) admits an AMLE relative to \(2^{X \setminus A}\).
\end{remark}

\begin{remark}
Let $X$ be a metric space, $A \subset X$ and let $f: A \to \mathbb R$ be a Lipschitz function. 
Fix $\lambda \in (0, 1)$ and define $m, M$ as in \eqref{mM}.
Then the convex combination
\[
    f_\lambda(x) := \lambda m(x) + (1 - \lambda) M(x), \quad \text{ for every } x \in X
\]
is a minimal Lipschitz extension of $f$ to $X$. 
However, note that not every function lying between $m$ and $M$ is a minimal Lipschitz extension.
For instance, let $X := \{ 0, \, 1, \, 3, \, 4 \} \subset \mathbb R$ and $A = \{ 0, \, 1 \}$.
Consider $f: A \to \mathbb R$ defined by
\[
f(0)=0 \quad \text{ and } \quad f(1)=1.
\]
Then \(\operatorname{Lip}(f,A)=1\). 
Recall that 
\[
m(x) = \sup_{a\in A} \{ f(a)-|x-a| \}\quad
 \text{ and } \quad M(x) = \inf_{a\in A} \{ f(a)+|x-a| \}.
\]
A direct computation shows that
\begin{align*}
& m(0) =0, \ m(1)=1,\ m(3)=-1 \, \text{ and } \, m(4)=-2,
\\
& M(0) =0,\ M(1)=1,\ M(3)=3  \, \text{ and } \, M(4)=4.
\end{align*}
Define \(g:X\to\mathbb R\) by
\[
g(0)=0,\quad g(1)=1,\quad g(3)=-1 \, \text{ and  } \, g(4)=4.
\]
It is straightforward to check that \(m(x)\le g(x)\le M(x)\) for every $x \in X$.
However
\[
\frac{|g(4)-g(3)|}{|4-3|}
=\frac{|4-(-1)|}{1}
=5>1=\operatorname{Lip}(f,A),
\]
so \(\operatorname{Lip}(g,X)>\operatorname{Lip}(f,A)\).
Therefore, although \(g\) lies between $m$ and $M$, it is not a minimal Lipschitz extension of $f$ to $X$.
\end{remark}

\begin{remark}\label{rmk01}
In contrast to the Euclidean setting (cf. e.g \cite{SC_2010, Jensen}), uniqueness of the discrete AMLE fails. 
The following example illustrates this. 
Let us consider $X = \{ x_1, \, x_2, \, x_3, \, x_4 \}$ equipped with the metric given as in the following graph.

\medskip

\begin{figure}[H]\centering
\begin{tikzpicture}[scale=1.3, every node/.style={inner sep=0pt}]
  % Nodes (horizontal)
  \node (a1) at (0,0) {$\bullet$};
  \node (a2) at (2,0) {$\bullet$};
  \node (a3) at (4,0) {$\bullet$};
  \node (a4) at (6,0) {$\bullet$};

  % Labels under each node
  \node at (0,-0.4) {$x_1$};
  \node at (2,-0.4) {$x_2$};
  \node at (4,-0.4) {$x_3$};
  \node at (6,-0.4) {$x_4$};

  % Curved edges from below
  \draw[-] (a1) to[bend right=20] node[midway, below, yshift=-4pt] {\small 1} (a2);
  \draw[-] (a2) to[bend right=20] node[midway, below, yshift=-4pt] {\small 2} (a3);
  \draw[-] (a3) to[bend right=20] node[midway, below, yshift=-4pt] {\small 1} (a4);

  % Curved edges from above
  \draw[-] (a1) to[bend left=20] node[midway, above, yshift=4pt] {\small 3} (a3);
  \draw[-] (a2) to[bend left=20] node[midway, above, yshift=4pt] {\small 3} (a4);

  % Wide curve below
  \draw[-] (a1) to[bend right=40] node[midway, below, yshift=-4pt] {\small 4} (a4);

\end{tikzpicture}
\caption{An example of a graph with four nodes.}
\label{fig:example-graph}
\end{figure}
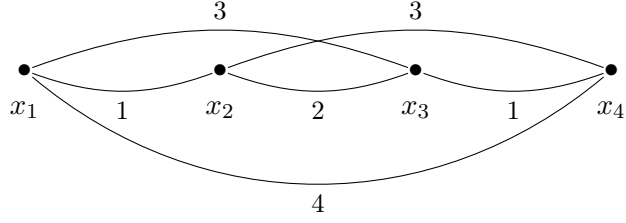
Consider $A = \{ x_1, \, x_2 \}$ and let $f: A \to \mathbb R$ be defined by $f(x_1) = 0$ and $f(x_2) = 1$ (and so $\Lip(f, A) = 1$).
For any fixed $t \in [-1, 3]$, define $h_t : X \to \mathbb R$ 
\begin{align*}
	h_t(x_1) = 0, \quad h_t(x_2) = 1 \quad \text{ and } \quad h_t(x_3) = h_t(x_4) = t.
\end{align*}
A direct computation shows that
\begin{align*}
	\dfrac{|h_t(x_i) - h_t(x_j)|}{d(x_i, x_j)} \leq 1, \quad \text{ for every } i, j \in \{1, 2, 3, 4 \},
\end{align*}
and therefore  each $h_t$ is a minimal Lipschitz extension of $f$.
Further, since $\Lip(h_t, \{ x_3, \, x_4 \}) = 0$, each $h_t$ is an AMLE of $f$ relative to $2^{X \setminus A}$.
This shows that uniqueness of AMLE fails in the discrete setting.
\end{remark}

\begin{remark}\normalfont
In Euclidean spaces (more generally in metric length spaces), the value function of $\epsilon$--tug--of-war games converges to AMLEs as $\epsilon \to 0^+$ and with a suitable metric, the two notions coincide, see \cite{BCR_2023, Gruyer, MRT_2012, PSSW_2009}.
For general graphs with finite width, this value function satisfies the so--called discrete infinity Laplace equation, see \cite{HW_2025}.
In contrast to the classical case, we will see that such value function may not even be a minimal Lipschitz extension.

\medskip

Let us recall the definition of the discrete infinity Laplacian on graphs. 
Let $G = (V, E)$ be a graph and let $A \subset V$.
We write $x \sim y$ if there is an edge connecting $x$ and $y$.
For given data \(f:A\to\mathbb{R}\), a function \(u:X\to\mathbb{R}\) satisfying
\begin{equation}\label{dis.infty}
\begin{cases}
u(x) = \dfrac{1}{2} \left(\displaystyle\max_{y \sim x} u(y)
+\displaystyle\min_{y \sim x} u(y) \right), & x\in X \setminus A,\\[2mm]
u(x)=f(x), & x\in A,
\end{cases}
\end{equation}
is a solution of the \textit{discrete infinity Laplace equation}.
According to \cite[Theorem 1.1]{HW_2025}, there exists a unique solution to the above functional equation.

\medskip

Coming back to our problem, let $X$ be a finite metric space. Consider the graph $G = (V, E)$, where $V = X$ and $E = X \times X$.
One can directly check that the solution to \eqref{dis.infty} is 
\begin{align*}
    u(x) := 
    \begin{dcases}
        \qquad\qquad f(x), & x \in A, \\
        \dfrac{1}{2} \left( \max_{a \in A} f(a) + \min_{a \in A} f(a) \right), & x \in X \setminus A
    \end{dcases}.
\end{align*}
This solution is in fact not necessarily a minimal Lipschitz extension of $f$ to the entire space. 
In particular, let $X=\{x_1, \, x_2, \, x_3, \, x_4\}$ be as in Remark~\ref{rmk01} and set $A=\{x_1, \, x_3\}$.
Consider $f:A\to\mathbb{R}$ with $f(x_1)=0$ and $f(x_3)=1$.
Then the solution $u:X\to\mathbb{R}$ of \eqref{dis.infty} will be $u=f$ on $A$ and
$u(x_2)=u(x_4)=\tfrac12$.
Clearly
\[
\Lip(f,A)=\frac{1}{3}
\qquad\text{while}\qquad
\Lip(u,X)=\frac{1}{2}>\Lip(f,A),
\]
so $u$ is not a minimal Lipschitz extension.
\end{remark}

\subsection{Proof of Theorem \ref{thm0}}
Consider $X = \{x_1, \, x_2, \, x_3, \, x_4, \, x_5\}$ and  $A = \{x_2, \, x_4\}$.  
The metric $d$ is given in Table~\ref{tab:metric-example2}.
\begin{table}[H]
\centering
\renewcommand{\arraystretch}{1.5} % <-- increases row height
\[
\begin{array}{c|ccccc}
	d(\cdot,\cdot) & x_1 & x_2 & x_3 & x_4 & x_5 \\
	\hline
	x_1 & 0   & 2 & 3 & 5 & 6 \\
	x_2 & 2 & 0   & 1 & 3 & 4 \\
	x_3 & 3 &  1 & 0   & 2 & 3 \\
	x_4 & 5 & 3 & 2 & 0   & 1 \\
	x_5 & 6 & 4 & 3 & 1 & 0
\end{array}
\]
\caption{Distance matrix}
\label{tab:metric-example2}
\end{table}

Let $f: A \to\mathbb R$ be given by $f(x_2) = 0$ and $f(x_4) = 1$. 
Thus $\Lip(f, A)= |0 - 1|/d(x_2, x_4) = 1/3$.
Recall that the McShane--Whitney extensions $m$ and $M$ in \eqref{mM}  are respectively the (pointwise) smallest and largest minimal Lipschitz extensions of $f$.
A direct computation yields
\begin{align*}
    & m(x_1) = - \dfrac{2}{3}, \quad m(x_3) = \dfrac{1}{3} \quad \text{ and } \quad m(x_5) = \dfrac{2}{3}, \\
    & M(x_1) = \dfrac{2}{3}, \quad M(x_3) = \dfrac{1}{3} \quad \text{ and } \quad M(x_5) = \dfrac{4}{3}.
\end{align*}
Therefore, any minimal Lipschitz extension $h: X \to \mathbb{R}$ of $f$ must satisfy $\Lip(h, X) = 1/3$ and
		\begin{align*}
				h(x_1) \in \left[ - \dfrac{2}{3}, \, \frac{2}{3} \right],\quad 	h(x_3) = \dfrac{1}{3} \quad \text{and}\quad 	h(x_5) \in \left[ \dfrac{2}{3}, \, \dfrac{4}{3} \right].
		\end{align*}
%From this, one can easily verify that 
%	\[
%		 |h(x_3) - h(x_4)| \leq 0.7 < 1.6, \quad
%		|h(x_4) - h(x_5)| \leq 0.8 < 1.6 \quad \text{and}\quad
%		|h(x_3) - h(x_5)| \leq 1.3 < 2.
%	\]

\medskip

Assume that $h$ is an AMLE relative to $2^{X \setminus A}$.
We consider two cases and show that both lead to a contradiction.

\medskip

\textbf{Case 1:} $\mathit{h(x_1) \neq  h(x_5)}$.  
Consider $U_{15} := \{x_1, \, x_5\} \subseteq X \setminus A$ and define $g: X \to \mathbb{R}$ by
\[
g(x) = h(x) \quad \text{ for every } x\in\{x_2, \, x_3, \, x_4\} \quad \text{ and } \quad
g(x_1) = g(x_5) = \dfrac{2}{3}. 
\]
Notice that $\Lip(g, U_{15}) = 0$ and $\Lip(g, X \setminus U_{15}) = \Lip(h, X \setminus U_{15}) \leq 1/3$.
Moreover, a direct computation gives:
\begin{alignat*}{3}
    & \dfrac{|g(x_1) - g(x_2)|}{d(x_1, x_2)} = \dfrac{1}{3}, & \quad\qquad & \dfrac{|g(x_5) - g(x_2)|}{d(x_5, x_2)} = \dfrac{1}{6} \leq \dfrac{1}{3},  \\
    & \dfrac{|g(x_1) - g(x_3)|}{d(x_1, x_3)} = \dfrac{1}{9} \leq \dfrac{1}{3}, & \quad\qquad & \dfrac{|g(x_5) - g(x_3)|}{d(x_5, x_3)} = \dfrac{1}{9} \leq \dfrac{1}{3}, \\
     & \dfrac{|g(x_1) - g(x_4)|}{d(x_1, x_4)} = \dfrac{1}{15} \leq \dfrac{1}{3}, & \quad\qquad & \dfrac{|g(x_5) - g(x_4)|}{d(x_5, x_4)} = \dfrac{1}{3}.\\
\end{alignat*}
Therefore $\Lip(g, X) = 1/3$ and so we deduce that $g$ is a minimal Lipschitz extension of $f$ and $g \big\vert_{X \setminus U_{15}} = h \big\vert_{X\setminus U_{15}}$.
However, since $h(x_1) \neq h(x_5)$, we have
\[
\Lip(g, U_{15}) = 0 
< \frac{|h(x_1) - h(x_5)|}{d(x_1, x_5)} 
= \Lip(h, U_{15}),
\]
This contradicts condition \textit{(ii)} of Definition \ref{def.amle.juuti}.  

\medskip

\textbf{Case 2:} $\mathit{h(x_1) = h(x_5)}$.  
Since we have
\[
    h(x_1) \in \left[ - \dfrac{2}{3}, \, \dfrac{2}{3} \right] \quad \text{ and } h(x_5) \in \left[  \dfrac{2}{3}, \, \dfrac{4}{3} \right],
\]
the equality forces $h(x_1) = h(x_5) = 2/3$.
Consider $U_{13} := \{ x_1, \, x_3 \}$ and  define $g: X \to \mathbb{R}$ by
\[
g(x) = h(x) \quad \text{ for every } x\in\{x_2, \, x_4, \, x_5\}, \quad  g(x_1) = \dfrac{1}{2}
\quad \text{ and  } \quad g(x_3) = \dfrac{1}{3}. 
\]
Note that $\Lip(g, U_{13}) = 1/18 \leq 1/3$ and $\Lip(g, X \setminus U_{13}) = \Lip(h, X \setminus U_{13}) \leq 1/3$.
Moreover, direct computation yields:
\begin{alignat*}{3}
    & \dfrac{|g(x_1) - g(x_2)|}{d(x_1, x_2)} = \dfrac{1}{4} \leq \dfrac{1}{3}, & \quad\qquad & \dfrac{|g(x_3) - g(x_2)|}{d(x_3, x_2)} = \dfrac{1}{3},  \\
    & \dfrac{|g(x_1) - g(x_4)|}{d(x_1, x_4)} = \dfrac{1}{10} \leq \dfrac{1}{3}, & \quad\qquad & \dfrac{|g(x_3) - g(x_4)|}{d(x_3, x_4)} = \dfrac{1}{3}, \\
     & \dfrac{|g(x_1) - g(x_5)|}{d(x_1, x_5)} = \dfrac{1}{36} \leq \dfrac{1}{3}, & \quad\qquad & \dfrac{|g(x_3) - g(x_5)|}{d(x_3, x_5)} = \dfrac{1}{9} \leq \dfrac{1}{3}.\\
\end{alignat*}
Therefore $\Lip(g, X) = 1/3$ and so we deduce again $g$ is a minimal Lipschitz extension of $f$ and  $g \big\vert_{X \setminus U_{13}} = h \big\vert_{X\setminus U_{13}}$.
However, we have
\[
\Lip(g, U_{13}) = \dfrac{1}{18}
< \dfrac{1}{9}
= \Lip(h, U_{13}),
\]
which contradicts condition \textit{(ii)} of Definition \ref{def.amle.juuti}.
\smallbreak 
Therefore no such AMLE relative to $2^{X \setminus A}$ can exist. \qed

\begin{remark}\label{rdms} $ \, $ \\
% \textit{(i)}   In the proof of Theorem \ref{thm0}, the constructed metric space embeds into $\mathbb R^d$ for every $d \geq 1$.
One can also construct a counterexample in a non-Euclidean setting.
This is illustrated in the following example.

\medskip

\noindent Consider $X = \{x_1, \, x_2, \, x_3, \, x_4, \, x_5\}$ and  $A = \{x_1, \, x_2\}$.  
The metric $d$ is given in Table~\ref{tab:metric-example}.
\begin{table}[H]
\centering
\renewcommand{\arraystretch}{1.5} % <-- increases row height
\[
\begin{array}{c|ccccc}
	d(\cdot,\cdot) & x_1 & x_2 & x_3 & x_4 & x_5 \\
	\hline
	x_1 & 0   & 2.0 & 0.6 & 1.1 & 1.7 \\
	x_2 & 2.0 & 0   & 1.6 & 1.1 & 0.5 \\
	x_3 & 0.6 & 1.6 & 0   & 1.6 & 2.0 \\
	x_4 & 1.1 & 1.1 & 1.6 & 0   & 1.6 \\
	x_5 & 1.7 & 0.5 & 2.0 & 1.6 & 0
\end{array}
\]
\caption{Distance matrix}
\label{tab:metric-example}
\end{table}

\noindent Let $f: A \to\mathbb R$ be given by $f(x_1) = 0$, $f(x_2) = 2$. 
Thus  $\Lip(f, A)= 1$.
Due to the McShane--Whitney bounds, any minimal Lipschitz extension $h: X \to \mathbb{R}$ must satisfy
		\begin{align*}
				h(x_3) \in [0.4, 0.6],\quad 	h(x_4) \in [0.9, 1.1] \quad \text{and}\quad 	h(x_5) \in  [1.5, 1.7].
		\end{align*}
From this, one can easily obtain that
	\[
		 |h(x_3) - h(x_4)| \leq 0.7 < 1.6, \quad
		|h(x_4) - h(x_5)| \leq 0.8 < 1.6 \quad \text{and}\quad
		|h(x_3) - h(x_5)| \leq 1.3 < 2.
	\]

\medskip

Assume that $h$ is an AMLE relative to $2^{X \setminus A}$.
We consider two cases according to the value of $h(x_4)$ and show that both lead to contradictions.  

\medskip

\textbf{Case 1:} $\mathit{h(x_4)\in (0.9,1.1] }$.  
Consider $U_{34} = \{x_3, \, x_4\} \subseteq X \setminus A$ and define $g: X \to \mathbb{R}$ by
\[
g(x) = h(x) \quad \text{ for every } x\in\{x_1, \, x_2, \, x_3, \, x_5\} \quad \text{ and } \quad
g(x_4) = 0.9. 
\]
It is straightforward to check that $g$ is $1$-Lipschitz.
Moreover, we have
\[
\Lip(g, U_{34}) = \frac{|h(x_3) - 0.9|}{1.6} 
< \frac{|h(x_3) - h(x_4)|}{1.6} 
= \Lip(h, U_{34}),
\]
since $h(x_4)>0.9$. 
This contradicts condition \textit{(ii)} of Definition \ref{def.amle.juuti}.  

\medskip

\textbf{Case 2:} $\mathit{h(x_4)=0.9}$.  
Consider $U_{45} = \{x_4, \, x_5\} \subseteq X \setminus A$ and define $g: X \to \mathbb{R}$ by
\[
g(x) = h(x) \quad \text{ for every } x\in\{x_1, \, x_2, \, x_3, \, x_5\}, 
\quad \text{ and  } \quad g(x_4) = 1.1. 
\]
Again one can check that $g$ is $1$-Lipschitz.
Furthermore,
\[
\Lip(g, U_{45}) = \frac{|1.1 - h(x_5)|}{1.6} 
< \frac{|0.9 - h(x_5)|}{1.6} 
= \Lip(h, U_{45}),
\]
since $h(x_5)\in[1.5,1.7]$. This also contradicts condition \textit{(ii)} of Definition \ref{def.amle.juuti}.  
\smallbreak 
Therefore no such AMLE relative to $2^{X \setminus A}$ can exist.
\end{remark}

\nocite{*}
\bibliographystyle{siam}
\bibliography{references}

\begin{thebibliography}{10}

\bibitem{SC_2010}
{\sc S.~N. Armstrong and C.~K. Smart}, {\em An easy proof of {J}ensen's theorem
  on the uniqueness of infinity harmonic functions}, Calc. Var. Partial
  Differential Equations, 37 (2010), pp.~381--384.

\bibitem{A_65}
{\sc G.~Aronsson}, {\em Minimization problems for the functional {${\rm
  sup}\sb{x}\,F(x,\,f(x),\,f\sp{\prime} (x))$}}, Ark. Mat., 6 (1965),
  pp.~33--53.

\bibitem{A_66}
\leavevmode\vrule height 2pt depth -1.6pt width 23pt, {\em Minimization
  problems for the functional {${\rm sup}\sb{x}\, F(x, f(x),f\sp\prime (x))$}.
  {II}}, Ark. Mat., 6 (1966), pp.~409--431.

\bibitem{A_67}
\leavevmode\vrule height 2pt depth -1.6pt width 23pt, {\em Extension of
  functions satisfying {L}ipschitz conditions}, Ark. Mat., 6 (1967),
  pp.~551--561.

\bibitem{ACJ_2004}
{\sc G.~Aronsson, M.~G. Crandall, and P.~Juutinen}, {\em A tour of the theory
  of absolutely minimizing functions}, Bull. Amer. Math. Soc. (N.S.), 41
  (2004), pp.~439--505.

\bibitem{BCR_2023}
{\sc L.~Bungert, J.~Calder, and T.~Roith}, {\em Uniform convergence rates for
  {L}ipschitz learning on graphs}, IMA J. Numer. Anal., 43 (2023),
  pp.~2445--2495.

\bibitem{C_2019}
{\sc J.~Calder}, {\em Consistency of {L}ipschitz learning with infinite
  unlabeled data and finite labeled data}, SIAM J. Math. Data Sci., 1 (2019),
  pp.~780--812.

\bibitem{CP_2007}
{\sc T.~Champion and L.~De~Pascale}, {\em Principles of comparison with
  distance functions for absolute minimizers}, J. Convex Anal., 14 (2007),
  pp.~515--541.

\bibitem{CEG_2001}
{\sc M.~G. Crandall, L.~C. Evans, and R.~F. Gariepy}, {\em Optimal {L}ipschitz
  extensions and the infinity {L}aplacian}, Calc. Var. Partial Differential
  Equations, 13 (2001), pp.~123--139.

\bibitem{DLV_2024}
{\sc A.~Daniilidis, M.~T. L{\^e}, and F.~Venegas}, {\em Absolutely minimal
  semi--{L}ipschitz extensions}, Calc. Var. Partial Differential Equations, 64
  (2025).

\bibitem{DJS_2019}
{\sc E.~Durand-Cartagena, J.~A. Jaramillo, and N.~Shanmugalingam}, {\em
  Existence and uniqueness of {$\infty$}-harmonic functions under assumption of
  {$\infty$}-{P}oincar\'e{} inequality}, Math. Ann., 374 (2019), pp.~881--906.

\bibitem{HW_2025}
{\sc F.~Han and T.~Wang}, {\em Discrete infinity {L}aplace equations on graphs
  and tug-of-war games}, Calc. Var. Partial Differential Equations, 64 (2025),
  p.~Paper No. 40.

\bibitem{Jensen}
{\sc R.~Jensen}, {\em Uniqueness of {L}ipschitz extensions: minimizing the sup
  norm of the gradient}, Arch. Rational Mech. Anal., 123 (1993), pp.~51--74.

\bibitem{J_2002}
{\sc P.~Juutinen}, {\em Absolutely minimizing {L}ipschitz extensions on a
  metric space}, Ann. Acad. Sci. Fenn. Math., 27 (2002), pp.~57--67.

\bibitem{JS_2006}
{\sc P.~Juutinen and N.~Shanmugalingam}, {\em Equivalence of {AMLE}, strong
  {AMLE}, and comparison with cones in metric measure spaces}, Math. Nachr.,
  279 (2006), pp.~1083--1098.

\bibitem{Gruyer}
{\sc E.~Le~Gruyer}, {\em On absolutely minimizing {L}ipschitz extensions and
  {PDE} {$\Delta_\infty(u)=0$}}, NoDEA Nonlinear Differential Equations Appl.,
  14 (2007), pp.~29--55.

\bibitem{MPR_2010}
{\sc J.~J. Manfredi, M.~Parviainen, and J.~D. Rossi}, {\em An asymptotic mean
  value characterization for {$p$}-harmonic functions}, Proc. Amer. Math. Soc.,
  138 (2010), pp.~881--889.

\bibitem{MRT_2012}
{\sc J.~M. Maz\'on, J.~D. Rossi, and J.~Toledo}, {\em On the best {L}ipschitz
  extension problem for a discrete distance and the discrete
  {$\infty$}-{L}aplacian}, J. Math. Pures Appl. (9), 97 (2012), pp.~98--119.

\bibitem{McShane}
{\sc E.~J. McShane}, {\em Extension of range of functions}, Bull. Amer. Math.
  Soc., 40 (1934), pp.~837--842.

\bibitem{PSSW_2009}
{\sc Y.~Peres, O.~Schramm, S.~Sheffield, and D.~B. Wilson}, {\em Tug-of-war and
  the infinity {L}aplacian}, J. Amer. Math. Soc., 22 (2009), pp.~167--210.

\bibitem{Whitney}
{\sc H.~Whitney}, {\em Analytic extensions of differentiable functions defined
  in closed sets}, Trans. Amer. Math. Soc., 36 (1934), pp.~63--89.

\end{thebibliography}
\end{document}